\documentclass[11pt]{amsart} 
\usepackage{anysize}
\usepackage{amsfonts} 
\usepackage{amssymb}
\usepackage[all]{xy}
\usepackage{epsfig}
\usepackage{url}
\usepackage{color}

\newtheorem{thm}{Theorem}
\newtheorem{cor}[thm]{Corollary} 
\newtheorem{prop}[thm]{Proposition}
\newtheorem{lemma}[thm]{Lemma}
\theoremstyle{definition}
\newtheorem{definition}{Definition}

\newenvironment{fig}[1][]{\mbox{}\\[.5\baselineskip]
    \refstepcounter{figure}
    \label{#1}
    \begin{minipage}{120mm}
      \small \noindent
      \textbf{Figure~\thefigure: }}
    {\end{minipage}
    \vspace{\baselineskip}}
\newcommand{\Fig}[1]{\mbox{}\\[.5\baselineskip]
    \refstepcounter{figure}
    \label{#1}
    \small \noindent
    \textbf{Figure~\thefigure}}

\newcommand{\Omit}[1]{{}}
\newcommand{\R}{{\mathbb R}}
\newcommand{\Z}{{\mathbb Z}}
\newcommand{\sm}{\setminus}
\newcommand{\lra}{\longrightarrow}
\newcommand{\sse}{\subseteq}

\newcommand{\core}{\mbox{\rm core}}
\newcommand{\st}{\mbox{\rm star}}
\newcommand{\wdge}{\omega}
\newcommand{\trunc}{\tau}
\newcommand{\B}{\mathcal{B}}
\renewcommand{\star}{\mbox{\rm star}}
\newcommand{\pol}{{}^\Delta}
\newcommand{\ovv}{\overline V} 
\newcommand{\auxGraph}{G_\mathcal{B}^*}

\begin{document}
\title[Perles' conjecture]{Examples and Counterexamples for Perles' Conjecture}
\author[Haase]{Christian Haase}
\email{haase@math.duke.edu}
\address{Department of Mathematics, Duke University, 
Durham, NC 27708-0320, USA}
\author[Ziegler]{G\"unter M. Ziegler$^*$}
\thanks{$^*$%
Supported by Deutsche Forschungsgemeinschaft (DFG)}
\stepcounter{footnote}
\email{ziegler@math.tu-berlin.de}
\address{MA 6-2, Institute of Mathematics, TU Berlin, D-10623 Berlin, Germany}
\setlength{\leftmargini}{7.5mm}

\begin{abstract}
The combinatorial structure of a $d$-dimensional simple convex
polytope -- as given, for example, by the set of the $(d-1)$-regular
subgraphs of the facets -- can be reconstructed from its abstract graph.
However, no polynomial/efficient algorithm is
known for this task, although a polynomially checkable certificate for 
the correct reconstruction exists. 

A much stronger certificate would be given by the following
characterization of the facet subgraphs, conjectured by M.~Perles:
``\emph{The facet subgraphs of a simple $d$-polytope are exactly all
  the $(d-1)$-regular, connected, induced, non-separating
  subgraphs}.'' 

We present non-trivial classes of 
examples for the validity of Perles' conjecture: 
In particular, it holds for the duals of cyclic polytopes, and for the
duals of stacked polytopes. 

On the other hand, we observe that for any $4$-dimensional
counterexample, the boundary of the (simplicial) dual polytope $P\pol$
contains a $2$-complex without a free edge, and without
$2$-dimensional homology. Examples of such complexes are known; we
use a modification of ``Bing's house'' (two walls removed) to
construct explicit $4$-dimensional counterexamples to Perles'
conjecture.
\end{abstract}

\maketitle

\section{Introduction}

If $P$ is a $d$-dimensional simple polytope, then
its graph $G=G(P)$ is a $d$-regular, $d$-connected graph.
If $F$ is any facet of~$P$, then $G_F=G(F)$
is a $(d-1)$-regular, induced, connected, non-separating
subgraph of $P$.
Whether these properties already characterize
the subgraphs of facets of a simple polytope
was asked by Perles \cite{Perles} a long time ago:

{\it Does every $(d-1)$-regular, induced, connected, non-separating
subgraph of the graph of a simple $d$-polytope correspond to
a facet of~$P$?}

This question is important in the context of reconstruction of polytopes. 
(See \cite{K2} for a general source on reconstruction of polytopes,
and in particular for Perles' unpublished work in this field.)
Simple polytopes can theoretically 
be reconstructed from their graphs -- 
this was originally conjectured by Perles, and first
proved by Blind \& Mani \cite{BM},
and particularly elegantly by Kalai~\cite{K1}.
However, although Kalai's proof can be implemented
to reconstruct polytopes of reasonable size \cite{A},
the reconstruction is far from easy or efficient 
(in the theoretical or practical sense).

Using ideas from Kalai's paper \cite{K1},
Joswig, Kaibel \& K\"orner \cite{KK} derived a certificate
for the reconstruction of a simple polytope that can indeed
be checked in polynomial time. They do, however, not (yet)
have a polynomial time algorithm to find/construct such a certificate.
If Perles' question had a positive answer, then the
reconstruction would be much easier, as was noted by Kalai \cite{K1},
and by Achatz \& Kleinschmidt~\cite{A}. 
In particular, Perles' conjecture proposes an
efficient criterion for recognizing the facet subgraphs of simple polytopes,
and thus an easy check whether a given list of subgraphs
is indeed the complete list of facet subgraphs.

Here we first establish that 
Perles' conjecture is true for non-trivial classes
of polytopes, including all the duals of stacked polytopes, and the 
duals of cyclic polytopes (Section~\ref{sec:positive}). 

However, our main result is that, perhaps surprisingly,
Perles' conjecture is \emph{not true} in general.
To prove this, we first identify certain topological obstructions
for Perles' conjecture (Section~\ref{sec:obstruction}).
We then take (known) $2$-dimensional
complexes that realize this obstruction --
the simplest one to handle being known as
``Bing's house with two rooms'' (with two walls removed).
We describe fairly general methods to produce counterexamples
from such $2$-complexes. A small, explicit counterexample derived this
way will be presented as an electronic geometry model 
(\url{www.eg-models.de}).

\section{Versions of Perles' Conjecture}

We will distinguish two slightly different versions of Perles' 
conjecture, where the first and stronger version is the one stated
in \cite{K1}, see also \cite[Problem 3.13*]{Z1}.

\begin{definition}
  A simple $d$-polytope $P$ \emph{satisfies Perles' conjecture} if
  every induced, connected, $(d-1)$-regular, non-separating subgraph
  of~$G(P)$ is the graph of a facet of~$P$.
\\
  A simple $d$-polytope $P$ \emph{weakly satisfies Perles' conjecture} 
  if every induced, $(d-1)$-connected,
  $(d-1)$-regular, non-separating subgraph of~$G(P)$ is the graph of
  a facet of~$P$.
\end{definition}

If $P$ is a simple polytope, then its polar dual $P\pol$ is
simplicial. The vertices and edges of $P$ correspond to the facets
and ridges of $P\pol$. The different notions in Perles'
conjecture for $P$ have immediate translations to the combinatorics of~$P\pol$. 
This yields the ``puzzle'' reformulation of~\cite{BM}.
\begin{itemize}
\item An induced subgraph $H$ of the graph of $P$ corresponds to a
  collection of facets and ridges of $P\pol$ such that, whenever two
  adjacent facets belong to the collection, then their common ridge
  also is in the collection. We will present this collection of facets
  and ridges in terms of the pure simplicial $(d-1)$-complex $\Gamma(H)$ that it
  generates. The dual graph of this ``pseudomanifold'' is the one that 
  Perles' conjecture in the original form refers to.
\item The induced subgraph $H$ is ($d-1$)-regular if and only if every
  ($d-1$)-simplex in $\Gamma(H)$ has exactly one free ($d-2$)-simplex.
  (A simplex is \emph{free} in a simplicial complex if it is contained
  in exactly one other/larger face of the complex.)
\item The induced subgraph $H$ is connected if and only if the corresponding
  subcomplex $\Gamma(H)$ is dually connected, that is, any two ($d-1$)-simplices 
  can be joined by a chain of successively adjacent facets and
  ridges. (This is often called ``strongly connected'' for pseudomanifolds.)
\item The induced subgraph does not separate its complement (graph
  theoretically) if and only if the corresponding subcomplex does not separate
  its complement (topologically): this is true since
  the topological complement retracts to the ``dual block complex'' \cite{Mu}
  whose $1$-skeleton is the graph of the complement.  
\item If $H$ is the graph of a facet, then the
  corresponding subcomplex $\Gamma(H)$ is a vertex star: the collection of all
  facets that contain a given vertex.
\end{itemize}

\begin{definition}
  A simplicial $d$-polytope \emph{satisfies Perles' conjecture} 
  if the vertex stars are the only pure $(d-1)$-dimensional, dually connected,
  non-separating subcomplexes of its boundary for which every maximal
  simplex has exactly one free face.
  
  A simplicial $d$-polytope \emph{weakly satisfies Perles' conjecture} 
  if the vertex stars are the only pure $(d-1)$-dimensional, 
  dually $(d-1)$-connected,
  non-separating subcomplexes of its boundary for which every maximal
  simplex has exactly one free face.
\end{definition}

\section{Positive results}\label{sec:positive}

Before embarking on the construction of counterexamples, 
and as an indication why these need to be rather complicated,
we will here demonstrate
that Perles' conjecture is not only \emph{quite plausible}, but
also that it is \emph{true} on large classes of examples:
it is certainly not an ``unreasonable'' conjecture.

We may safely assume that Perles himself verified his conjecture
for non-trivial classes of examples; the Diplomarbeit of
Stolletz \cite{St} also has a number of positive results.
Thus, Perles' conjecture is true, for example, for all $d$-dimensional 
polytopes with $d\le 3$ 
(even without a restriction to simple/simplicial polytopes),
for all $4$-dimensional simple polytopes with at most $8$ vertices
\cite[Sect.~5.2.1]{St},
for all $4$-dimensional product polytopes \cite[Kap.~2]{St},
for all $4$-dimensional cyclic polytopes \cite[Kap.~4]{St}, and
for all stacked polytopes \cite[Satz~3.1.3]{St}.

Indeed, the treatment of vertex truncations, and thus of 
stacked polytopes, presents no greater difficulties.

\begin{prop}[{\rm \cite[Lemma 3.1.2]{St}}]
Truncation (``cutting off a simple vertex'')
preserves the validity of Perles' conjecture:
A vertex-truncated simple $d$-polytope $P':=\trunc_v(P)$ satisfies
Perles' conjecture if and only if $P$ satisfies it.
\end{prop}

\Omit{\begin{proof}
Let $w_1,\ldots,w_d$ be the neighbors of~$v$, 
then $G'=G(P')$ is constructed from $G=G(P)$ by replacing the
edges $(w_i,v)$ by new edges $(w_i,v_i)$, and adding a new 
$d$-clique $K\subset G$ (complete subgraph) 
on the vertex set $\{v_1,\ldots,v_d\}$. Conversely, $G$ arises from~$G'$ by
contraction of the edges of~$K$.
\begin{center}
\input{trunc.pstex_t}\begin{fig}[fig:trunc]
\end{fig}\end{center}
Now if $H'\subset G'$ is any connected, non-separating
$(d-1)$-regular subgraph, then three
different situations are possible.
First, $H'$ may contain no vertex of $K$, then it describes also
a subgraph of~$G$; hence it corresponds to a facet of~$P$, 
and thus also of~$P'$.
Secondly, $H'$ may contain all the vertices of~$K$, then it coincides
with $K$, and hence corresponds to the ``new'' facet of~$P'$.
Finally, if $H'$ contains some, but not all the vertices
of~$K$, then $(d-1)$-regularity implies that the intersection
with $K$ is a $(d-1)$-clique. Contraction of $K$ (and thus of
the $(d-1)$-clique) yields a new, $(d-1)$-regular subgraph,
which by Perles' conjecture is the graph of a facet of~$P$.
The validity of Perles' conjecture for $P'$ thus follows from the
validity for~$P$ -- and conversely.
\end{proof}
}

\begin{cor}[{\rm Stolletz \cite[Satz~3.1.3]{St}}]
All stacked polytopes (dually: the multiple 
vertex-\allowbreak trunca\-tions of simplices) 
satisfy Perles' conjecture.
\end{cor}

We now develop a rather systematic method
to prove Perles' conjecture for classes of simple polytopes.
We demonstrate its use for products, for wedges, and in
particular for the duals of cyclic polytopes.
Here is the idea: Let $P$ be a simple $d$-polytope, $G=G(P)$
its graph, with vertex set~$V$. We are interested in $(d-1)$-regular
induced connected subgraphs $H$ that don't separate;
as induced subgraphs they are given by their vertex sets $V_H\subset V$.
The `Ansatz' is to concentrate on the vertex set 
$\ovv=V\sm V_H$ of the complement of $H$ in~$G$.
In the \emph{first step} we (try to) classify all those vertex 
sets $\ovv\subset V$ for which  
\begin{enumerate}
\item[(1)]the induced subgraph $G[\ovv]$ is connected, and
\item[(2)]$\ovv$ is triangle-closed and quadrilateral-closed:
if $\ovv$ contains ``all but one'' vertices of a triangle or
quadrilateral of~$G$, then it must contain all its vertices.
\end{enumerate}
The first condition is necessary since we want $H$ to be non-separating.
If the second condition is violated, then the ``missing vertex'' 
has two neighbors in~$\ovv$, and thus has degree at most $d-2$ in~$H$.

In the \emph{second step} we then identify those sets $\ovv$
for which every vertex not in~$\ovv$ has exactly one neighbor 
in~$\ovv$. Our hope is then to end up with only
the (complements of) facet subgraphs, plus the trivial examples
given by~$\ovv=V$.

\begin{prop}\label{prop:product}
If two simple polytopes $P_1$ and $P_2$ satisfy Perles' conjecture,
then so does their product~$P_1\times P_2$.
\end{prop}

\begin{proof}
Let $H\subset G(P_1\times P_2)=G_1\times G_2$ be a $(d_1+d_2-1)$-regular
subgraph that is connected and does not separate ($d=d_1+d_2$).
Now (Step~1) if the complement vertex set $\ovv$ satisfies (1) and (2),
then it is a product set, $\ovv=\ovv_1\times\ovv_2$.
\begin{center}
\input{prod.pstex_t}\begin{fig}[fig:prod]
Property (2) at work in a product graph
\end{fig}
\end{center}
Now (Step~2) if we can choose $v_i\in V_i\setminus \ovv_i$ for $i=1,2$,
then this yields a vertex $(v_1,v_2)\in H$ 
that in~$G_1\times G_2$ has graph-theoretic distance at least~$2$
from the complement set $\ovv_1\times \ovv_2$
(symbolized by the vertex $\bullet$ in our figure).
Thus this vertex of~$H$ has degree $d_1+d_2$, which is impossible.

Hence we get that the $(d_1+d_2-1)$-regular, induced, connected and 
non-separating subgraphs all have the form $H=G_1\times H_2$ or
$H=H_1\times G_2$, where $H_i\subset G_i$ is $(d_i-1)$-regular, 
induced, connected and 
non-separating. By Perles' conjecture, which is valid for $P_i$,
it follows that $H_i$ is the graph of a facet of~$P_i$,
and hence that $H$ is the graph of a facet of~$P_1\times P_2$.
%
\end{proof}

We refer to Klee \& Walkup \cite{KlWa} and to Holt \& Klee \cite{HK1} 
for the construction of
the wedge $\wdge_FP$ of a polytope $P$ over a facet $F$.
(It is what you think it ought to be.)

\begin{prop}\label{prop:wedge}
If $P$ is a simple $d$-polytope that satisfies Perles' conjecture,
and $F$ is a facet of~$P$, then $\wdge_F(P)$ satisfies Perles' conjecture
(and conversely).
\end{prop}

\begin{center}
\input{wedge.pstex_t}\Fig{fig:wedge}
\end{center}

\begin{proof}
[Sketch of proof]The vertex set of $P':=\wdge_F(P)$ may be decomposed as
$V_F\cup V_T\cup V_B$, into the vertices that lie in the 
face $F$ resp.\ not in the bottom facet resp.\ not in the top facet of the 
wedge. In Step~1 one then verifies that if $\ovv$ meets both
$V_B$ and $V_T$, then it contains always ``none or both'' from
a pair of corresponding top and bottom vertices.
\end{proof}

\Omit{\begin{proof}
The vertex set of $P':=\wdge_F(P)$ may be decomposed as
$V_F\cup V_T\cup V_B$, into the vertices that lie in the 
wedge facet resp.\ not in the bottom facet resp.\ not in the top facet. 

Now let $H$ be a $d$-regular subgraph of $G'=G(P')$,
and consider the vertex set of the complement, $\ovv:=V'\setminus H'$. 
The interesting case occurs when $\ovv$ contains
vertices both from~$V_B$ and from~$V_T$. 
Then one derives that $\ovv$ contains
always ``none or both'' of any top-bottom pair:
here we consider a $\ovv$-path from a vertex in $v_B$ to
a vertex in $v_T$, and ``on the way'' use the property
that if all but one vertex of a triangle or quadrilateral lie in 
the $\ovv$, then they all do.
(This is illustrated by the following figure, where ``$\bullet$''
denotes a vertex in~$H$, while ``$\circ$'' represents a vertex in~$\ovv$.) 
\begin{center}
\input{wedge.pstex_t}\Fig{fig:wedge}
\end{center}
Furthermore, every vertex in if $v\in V_F$ has exactly
one neighbor in~$V_B$ and one in~$V_T$, and these are
adjacent. Thus if $v\in V_F$ is in~$H'$, then
also its two neighbors in $V_B$ resp.\ $V_T$ lie in~$H'$.
{}From this one derives that the restriction of
$H'$ to $V_F\cup V_B$ is $d$-regular, connected, and non-separating,
and Perles' conjecture can be applied to it.
\end{proof}
}

\begin{lemma}\label{lemma:cyclic-wedge}
If $d>2$ is odd, then $C_d(n)\pol$ is combinatorially equivalent to 
$\wdge_F C_{d-1}(n-1)\pol$, for a facet $F\cong C_{d-2}(n-2)\pol$.
If $d>1$ is even, and $n=d+1$, then $C_d(n)\pol$ is a $d$-simplex; for $n=d+2$
it is a product of two $\frac d2$-simplices.
\end{lemma}

\begin{proof}
To be derived from Gale's evenness criterion; see \cite[p.~62]{Gr}
or~\cite[p.~14]{Z1}.
\end{proof}

\begin{thm}
The cyclic polytopes $C_d(n)$ satisfy Perles' conjecture.
\end{thm}

\begin{proof}We work with the duals, $P= C_d(n)\pol$.
By Propositions \ref{prop:product} and \ref{prop:wedge} 
plus Lemma~\ref{lemma:cyclic-wedge},
we need only treat the case where $d=2e$ is even, and $n>d+2$.

(i.)
The graph $G=G(C_d(n)\pol)$ has a simple combinatorial description,
via Gale's evenness criterion: Its vertices $v$ 
correspond to those subsets $S\in\binom{[n]}{d}$ which
split into a union of $e$ adjacent pairs modulo~$n$ (that is,
we identify the elements of $[n]=\{1,\ldots,n\}$  with $\Z_n$).
The splitting into adjacent pairs is always unique.
We shall call a ``block'' any non-empty union of adjacent pairs
that is contiguous, that is, without a gap.
Two vertices are adjacent if they differ in a single element,
that is, if one arises from the other by moving one block by
``one unit.'' This also provides a canonical orientation on each
edge: we put a directed edge $v_S$\,\includegraphics{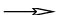}\,$v_T$ 
if we get from $S$ to $T$
by moving a block ``up'' (mod~$n$).
The resulting digraph is, of course, not acyclic.
\begin{center}\input{circle.pstex_t}\begin{fig}[fig:circle]
One single pair is moved, thus we get an edge of $G$,
and of $G'$, corresponding to $12356...\rightarrow13456...$.
\end{fig}
\end{center}
\vskip-3mm
We shall also consider the subgraph $G'\sse G$, which has the
same vertex set as~$G$, but only retains those directed edges
that correspond to moves of single pairs.
(In Figure~\ref{fig:moebius} below, this is the graph formed
by the straight edges only.)
The structure of the digraph $G'$ is closely linked to the poset
\[
L\ \ :=\ \ \{(j_1,j_2,\ldots,j_d)\in\Z^d:
j_{2k}   =  j_{2k-1}+1,\ 
j_{2k+1}\ge j_{2k}+1,\ 
j_d\le j_1+n-1\},
\]
equipped with componentwise partial order.
This poset is a distributive lattice.
Moreover, reduction modulo~$n$ defines a surjective,
and locally injective, digraph map
$\Phi:L\lra G'$, from the Hasse diagram of $L_d(n)$ onto the digraph~$G'$.

(ii.)
Now assume (Step~1) that $\ovv\subset V$ satisfies (1) and (2).
Every move of a block can be decomposed into a sequence of 
moves of pairs. Thus property (2) implies that also $G'[\ovv]$ is connected
(as an undirected graph): every directed arc in $G[\ovv]$ canonically
corresponds to a sequence of directed arcs in~$G'[\ovv]$.
Thus $G[\ovv]$ is acyclic if and only if $G'[\ovv]$ is acyclic.

(iii.) 
Next we treat the case that $G'[\ovv]$ contains a directed cycle.
Every such cycle lifts into a two-way infinite, maximal chain 
$C\sse \Phi^{-1}(\ovv)$ in the lattice~$L$. \emph{Every} 
element $w\in L$ is contained in a finite interval $[v',v'']$ 
between elements $v',v''\in C$. Now $C$ restricts to a maximal
chain $\gamma_0$ in the interval $[v',v'']$, while $w$ will lie on
some other maximal chain $\gamma$ of this interval. But in a distributive 
(and hence semimodular) lattice one can move from any maximal chain
to any other one by one-element exchanges, see Bj\"orner~\cite[Sect.~6]{Bj}.
This implies via (2) not only the elements of $\gamma_0$, but
also all elements of chains that we can move to, belong to $\Phi^{-1}(\ovv)$.
Thus we get $w\in\Phi^{-1}(\ovv)=L$, and hence $\ovv=V$.

\begin{center}
\input{moebius.pstex_t}\qquad\qquad\input{lattice.pstex_t}%
\begin{fig}[fig:moebius]
This depicts the digraph $G$, and the lattice~$L$, for $d=4$, $n=8$.
The digraph $G$ is finite, the ends are to be identified according
to the capital letters. It is not planar, but for $d=4$ it
embeds into a M\"obius band. The lattice~$L$ is infinite. 
\end{fig}
\end{center}

(iv.)
Thus we may assume that $G'[\ovv]$ is acyclic.
We can then identify $G[\ovv]$ with 
an induced subgraph of the Hasse diagram of~$L$: every connected component 
of $\Phi^{-1}(\ovv)$ is isomorphic to $G'[\ovv]$.
With this the properties (1) and (2), and lower and upper semi-modularity of
the lattice $L$, imply that $\ovv$ corresponds to an interval in~$L$
(see again \cite{Bj});
that is, there are elements $x_0,x_1\in L$
such that $G'[\ovv]$ is isomorphic to the Hasse diagram
of the interval $[x_0,x_1]\sse L$. This corresponds to
unique vertices $v_0,v_1\in\ovv$
such that $\ovv$ consists of all vertices of~$G'$
that lie on a directed path of minimal length from $v_0$ to $v_1$.

(v.)
Assume that $\ovv$ contains no ``no-gap vertex,'' 
whose set~$S$ consists of one single block of size~$d$. 
(In the case $d=4$, cf.\ Figure~\ref{fig:moebius},
this means that $\ovv$ contains no 
vertex on the border of the M\"obius strip.)
Then every no-gap vertex needs to have exactly one
neighbor in~$\ovv$, which is a one-gap vertex
consisting of $d$ elements from a block of $d+1$ adjacent
vertices. Every such block corresponds to a $(\frac{d}2+1)$-clique;
these cliques contain two no-gap vertices each, and each no-gap
vertex is contained in two of these $(\frac{d}2+1)$-cliques.
(In Figure~\ref{fig:moebius}, this corresponds to the chain
of triangles in the boundary of the M\"obius strip.)

Now we lift the situation to the lattice $L$, where the
chain of cliques gives rise to $\frac d2$ distinct chains,
which are disjoint by our assumption $n>d+2$.
The interval $[x_0,x_1]$ can contain at most one element
from each of the chains.
%
Thus it can contain at most
$\frac d2$ elements from the chains. Thus at most
$d$ no-gap vertices are adjacent to a vertex in~$[x_0,x_1]$.
Projecting this back to $G$, we find that at most
$d$ no-gap vertices are adjacent to a vertex in $\ovv$. But
they all have to be, so we get $n\le d$, a contradiction.

(vi.)
If $\ovv$ contains a one-block vertex, then 
we see, using (2), that both $v_0$ and $v_1$
must be one-block vertices. By symmetry, we may then
assume that $S_0=\{1,2,\ldots,d\}$ and $S_1=\{k+1,\ldots,d+k-1\}$
for some~$k$. This ends Step~1: and Step~2 
-- the identification of those parameters $k$ for which
each vertex in $V\sm\ovv$ has exactly one neighbor in~$\ovv$ -- 
is now easy.
\end{proof}

\section{The obstruction \label{sec:obstruction}}

We now restrict our discussion to the case $d=4$, where we will construct
a counterexample to the weak Perles conjecture.%
\renewcommand{\thefootnote}{*}%
\footnote{The pictures try to illustrate the situation in $d=3$, though.} 
{}From it, one gets counterexamples to the weak Perles conjecture
in all dimensions~$d>4$, from the constructions of Section~\ref{sec:positive},
such as wedges and products.

The discussion in this section motivates our construction; it 
leads us to well-characterized obstructions to the validity of 
Perles' conjecture (for any specific polytope). For simplicity,
we formulate this for $d=4$; the generalization to $d \ge 4$ is immediate.

\begin{prop}\label{prop:obstruction}
  Let $P^\Delta$ be a simplicial $4$-polytope whose boundary complex
  $\Delta:=\Delta(\partial P^\Delta)$ contains a dually connected pure
  $3$-dimensional subcomplex $\Gamma$, all whose tetrahedra have
  exactly one free triangle in $\Gamma$. 

  Then $\Gamma$ $($and in particular $\Delta)$ contains a pure
  $2$-dimensional, dually connected subcomplex $\core(\Gamma)$
  without free edges.
  Moreover, $\Gamma$ separates the boundary of 
  $P^\Delta$ if and only if $\core(\Gamma)$ does, and
  $\core(\Gamma)$ is empty if and only if $\Gamma$ is a vertex star.
\end{prop}

In the case $d=3$, $\core(\Gamma)$ would be $1$-dimensional,
with neither a cycle (non-separating) nor a free face: a leafless
tree. This (re)proves Perles' conjecture for $3$-polytopes.

\begin{proof}
  We specify a pure $2$-dimensional subcomplex $\core(\Gamma)$ of
  $\Gamma$ which will comply with our conditions (the definition is
  due to Carsten Schultz):
\emph{Every tetrahedron $\sigma \in \Gamma$ has a unique vertex
  $v(\sigma)$ opposite to its free face. A triangle $\tau$ belongs to
  $\core(\Gamma)$ if both its neighboring tetrahedra $\sigma_1,
  \sigma_2$ belong to $\Gamma$, and if $v(\sigma_1) \neq v(\sigma_2)$.}
  \begin{center}
    \includegraphics[height=35mm]{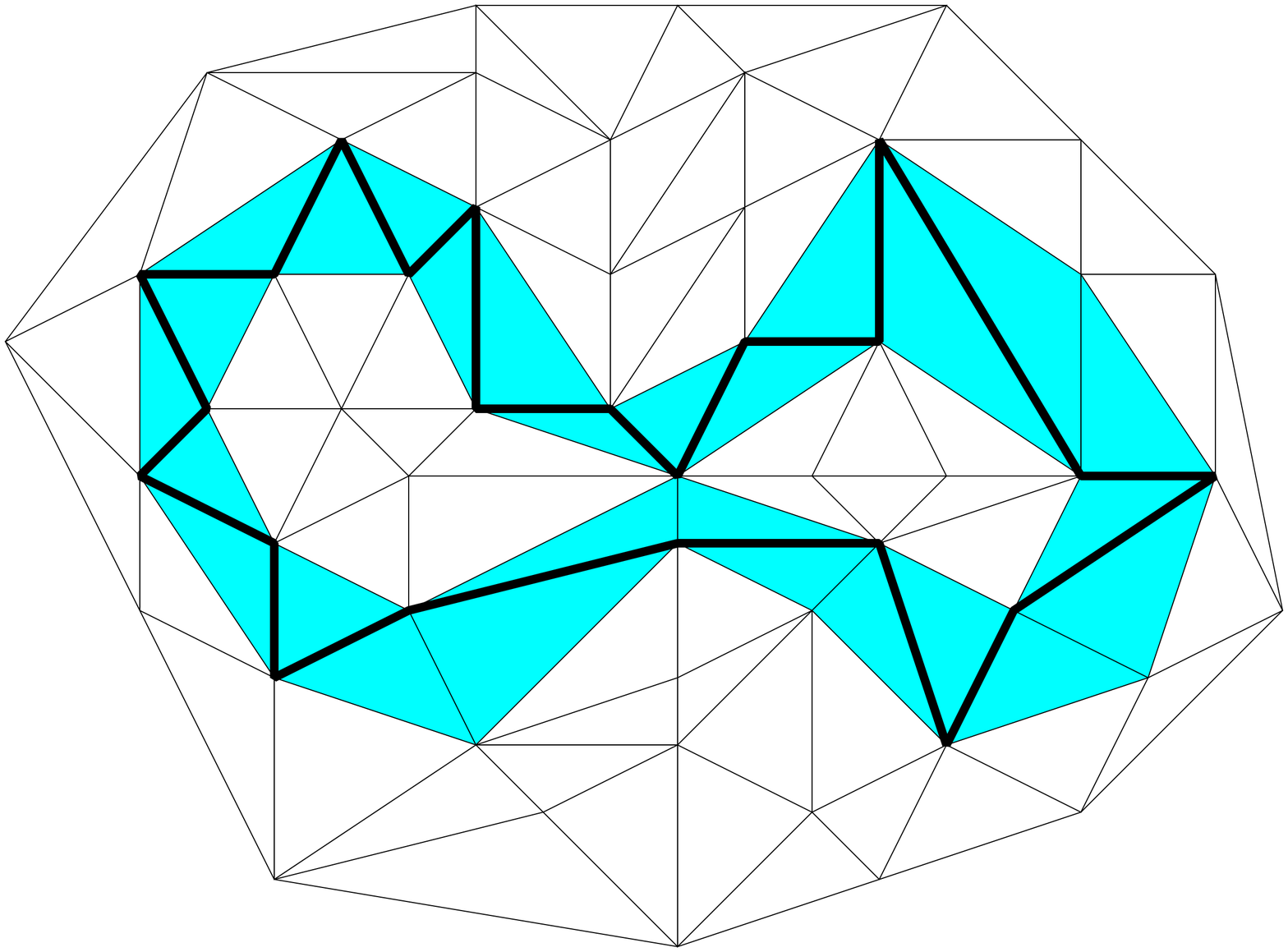}
    \begin{fig}[fig:core]
      If the grey $2$-complex is $\Gamma$, 
      then the black $1$-dimensional subcomplex is its core.
    \end{fig}
  \end{center}

\noindent{\em Claim 1. }
    The complex $\core(\Gamma)$ has no free edge. 
\smallskip

  Suppose that $\varrho$ is a free edge of the triangle $\tau \in
  \core(\Gamma)$. The tetrahedra of $\Delta$ which contain~$\varrho$
are cyclically ordered: $\sigma_1, \ldots, \sigma_n$, as
  in Figure~\ref{fig:cycle}.
  \begin{center}
    \input{labelledstar.pstex_t}
    \Fig{fig:cycle}
  \end{center}
  We can assume that $\tau = \sigma_1 \cap \sigma_2$ and $\sigma_3 \in
  \Gamma$. For, if $\sigma_3$ and $\sigma_n$ were both not in
  $\Gamma$, then this would imply $v(\sigma_1) = v(\sigma_2)$, and
  hence $\tau \not\in \core(\Gamma)$.

  Because $\sigma_2 \cap \sigma_3 \not\in \core(\Gamma)$, we have
  $v(\sigma_3) = v(\sigma_2)$. This vertex lies in
  $\sigma_2\cap\sigma_3$, and it cannot be the vertex of $\sigma_2$
  that is opposite to~$\tau$; thus we get that 
  $v(\sigma_3) = v(\sigma_2)\in \varrho$. 
  In particular, this means that $\sigma_3\cap\sigma_4$ is not the
  free face of~$\sigma_3$, and thus $\sigma_4 \in \Gamma$. 
  Iterating these arguments, we see that $\sigma_1,\ldots,\sigma_n \in
  \Gamma$ and $v(\sigma_n) = \ldots = v(\sigma_2) \neq v(\sigma_1)$.
  So, $\sigma_n \cap \sigma_1$ is another triangle in $\core(\Gamma)$.
\medskip

\noindent{\em Claim 2. }
    The complex 
    $\core(\Gamma)$ is empty if and only if $\Gamma$ is a vertex star.
\smallskip

  If $\core(\Gamma)$ is empty, as $\Gamma$ is dually connected, all
  the vertices $v(\sigma)$ are equal, say, $v_0$. Hence $\Gamma$ is
  part of $\star_\Delta(v_0)$. But then $\Gamma$ must be the whole
  star, because there is only one free face per tetrahedron.
\medskip

\noindent{\em Claim 3. }
    If $\core(\Gamma)$ separates then so does $\Gamma$.
\smallskip

  Alexander duality implies that a subcomplex of a $3$-sphere
  separates if and only if it has non-trivial $2$-dimensional
  homology.
  The $3$-dimensional complex $\Gamma$ collapses down to a
  $2$-dimensional subcomplex $\Gamma'$, which contains
  $\core(\Gamma)$. Since $\Gamma'$ has dimension $2$, it is clear that
  the (relative) homology group $H_3(\Gamma',\core(\Gamma))$ vanishes. 
  Hence, the map $\iota$ in the exact sequence
  \begin{equation*}
    \def\labelstyle{\displaystyle}
    \xymatrix@R=.7\baselineskip{
      H_3(\Gamma,\core(\Gamma)) \ar[r] \ar@{=}[d] &
      H_2(\core(\Gamma)) \ar[r]^-\iota & H_2(\Gamma) \ar[r] & \cdots
      \\
      H_3(\Gamma',\core(\Gamma)) \ar @{}[r]|{=0} &
      }
  \end{equation*}
  is injective, and the claim is proved.
\end{proof}

Let us further analyze the above situation. The aim is to reveal the
combinatorial properties of complexes that appear as a
core (in addition to the topological ones that we have
already seen). In Section~\ref{sec:example} we will show how to
guarantee these properties using stellar subdivisions.

Consider now the structure of $\Gamma$ locally, restricted
to the (open) star of a vertex $v_0$. This vertex star is
divided by~$\core(\Gamma)$ into several pieces and $\Gamma$ is a union
of some of these pieces:
\begin{lemma} \label{prop:pieces}
  Let $\sigma,\sigma' \in \star_\Delta(v_0)$ be tetrahedra in
  the same piece, that is, such that 
  there is a dual path in $\star_\Delta(v_0)$
  between $\sigma$ and $\sigma'$ that does not meet $\core(\Gamma)$.
  If $\sigma \in \Gamma$ and $v(\sigma)=v_0$, then $\sigma' \in \Gamma$
  and $v(\sigma')=v_0$.
\end{lemma}
\begin{cor} \label{cor:pieces}
  For every vertex $v_0$ of $\core(\Gamma)$ the complex which is
  generated by $\{ \sigma \in \Gamma : v(\sigma) = v_0 \}$ is
  the closure of a union of components of
  $\star_\Delta(v_0)\sm\core(\Gamma)$.

  Furthermore, two such components that correspond to different
  $v_0$'s intersect at most in codimension
  one, and if so, then this intersection is included in
  $\core(\Gamma)$.
\end{cor}
\begin{center}
  \includegraphics[width=60mm]{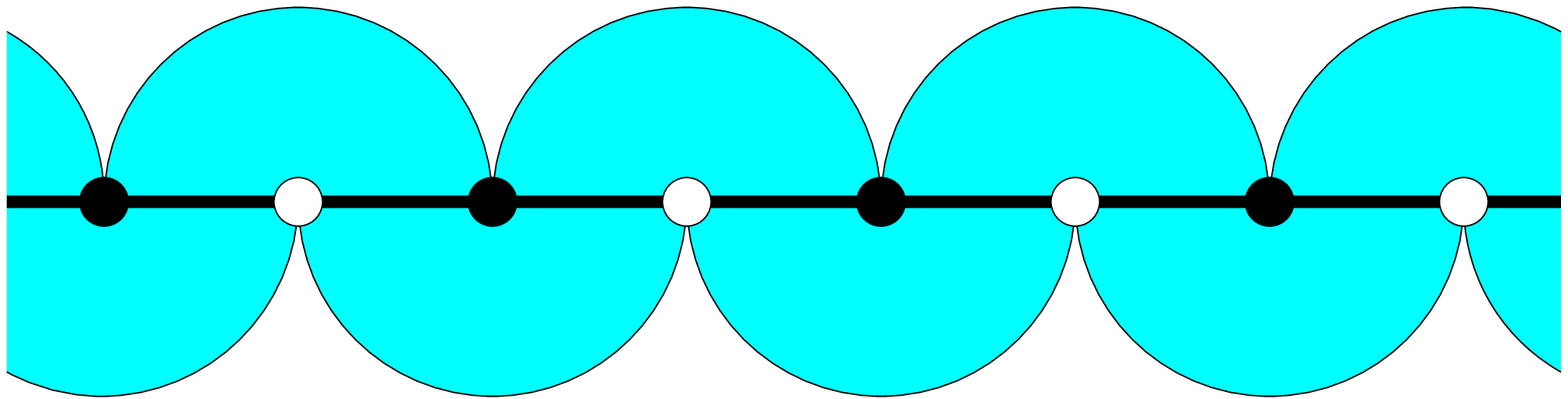}
\begin{fig}[fig:core2] illustrates Corollary~\ref{cor:pieces}\end{fig}
\end{center}

\section{A counterexample \label{sec:example}}

Proposition~\ref{prop:obstruction} makes one
ask for a $2$-dimensional simplicial complex $\Gamma$, 
without a free edge (every edge is contained in at least two triangles),
which does not separate (that is, $H_2(\Gamma,\Z)=\{0\}$), and
which is embeddable in~$\R^3$ (in particular, $H_1(\Gamma,\Z)$ has 
no torsion).
Such $2$-dimensional complexes do exist. The most prominent examples
are probably Borsuk's ``dunce hat'' \cite{Borsuk,Zeeman}, and Bing's
``house with two rooms'' \cite{Bing} \cite[p.~4]{Hatcher}, 
which are even contractible.

To keep the constructions for the following simpler, we shall
put an extra condition on the complexes we look at: we want them
to be ``essentially manifolds,'' that is, we want every point to
have a closed neighborhood that is homeomorphic either to two or to
three triangles that are joined together at a common edge.
(Thus, we admit singular curves, along which a manifold branches
into three parts, but we do not admit singular points that may be
more complicated than the points on a singular curve.)

\begin{center}
  \input{bing3a.pstex_t}
  \begin{fig}[fig:bing]
    The complex $\B$ has two rooms: The downstairs room is connected to
    the outside via the chimney through the upstairs room, the upstairs
    room is connected to the outside via the chimney through the
    downstairs room.
  \end{fig}
\end{center}

In the following, we provide an as-concrete-as-possible
description of a specific counterexample of this type.
Further 
counterexamples of the same type may be obtained along the same lines.

Our point of departure is a modification $\B$ of Bing's house. Two
extra walls are missing that would usually be added to make the
interior of each room simply connected (cf.~Figure~\ref{fig:bing}).

This $\B$ can be embedded into a pile of $2 \times 3 \times 4$ cubes -- a
cubical complex. Triangulate the pile by the arrangement
of hyperplanes $x_i-x_j=k$ to get a simplicial
complex (cf.~Figure~\ref{fig:cube}). Add a cone over the boundary to
get the boundary of a simplicial $4$-polytope $Q\pol$ (any
subdivision by a hyperplane arrangement is ``coherent''/``regular'';
cf.~\cite[Lecture~5]{Z1}).
\smallskip

\begin{center}
  \includegraphics[width=30mm]{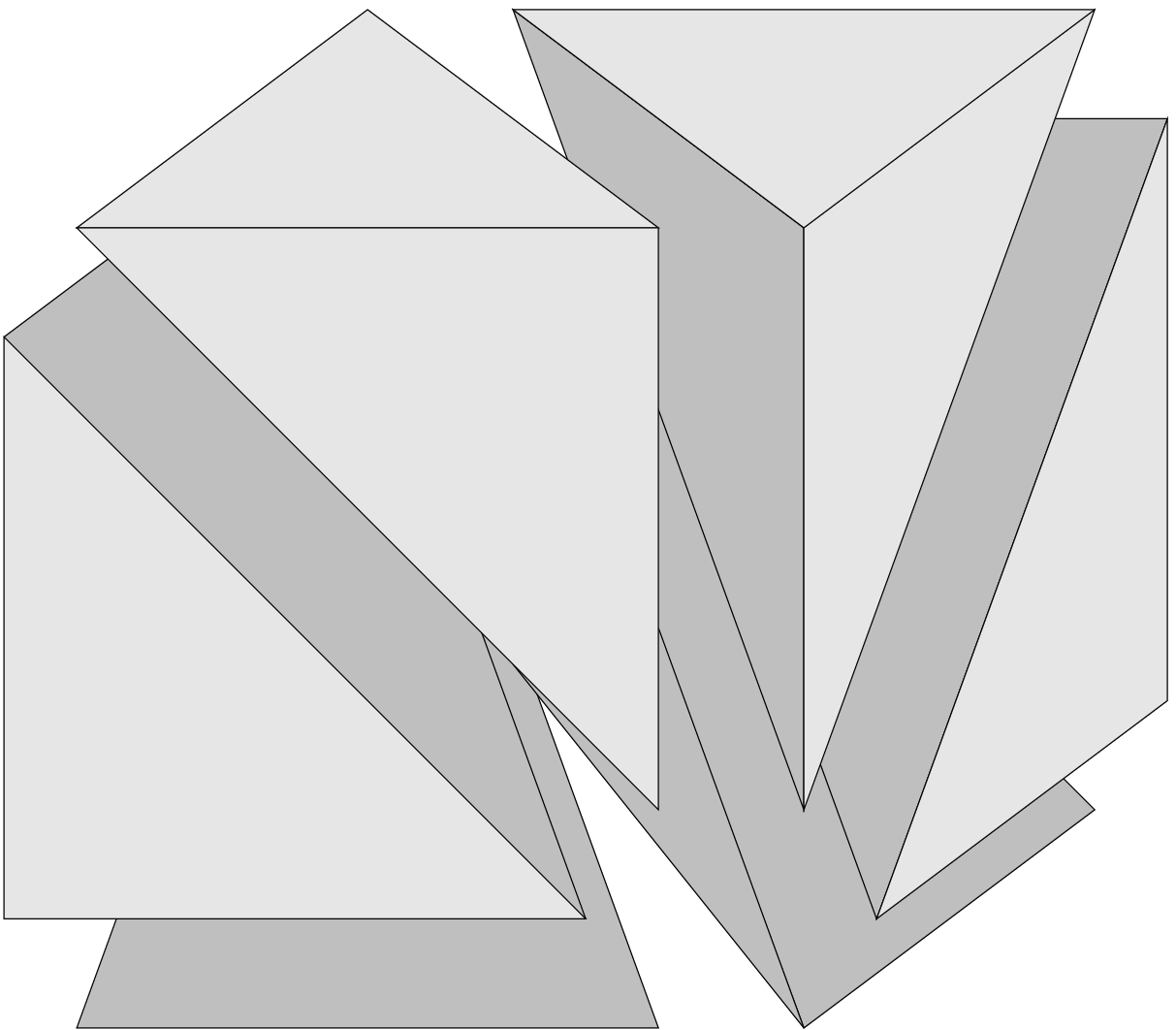}
  \Fig{fig:cube}
\end{center}
This polytope contains a triangulation of $\B$ as a subcomplex, which
we also denote by~$\B$. The vertices of this subcomplex 
now get ``colors'' $0$, $1$, or $2$, by assigning to each vertex the sum
of its coordinates modulo 3.
The only edges that join vertices of the same color are the
diagonals of the cubes, which are not in $\B$. The two chimneys
look as in Figure~\ref{fig:chimneys}.
\smallskip

\begin{center}
  \raisebox{9mm}{The upper chimney:}
  \includegraphics[width=25mm]{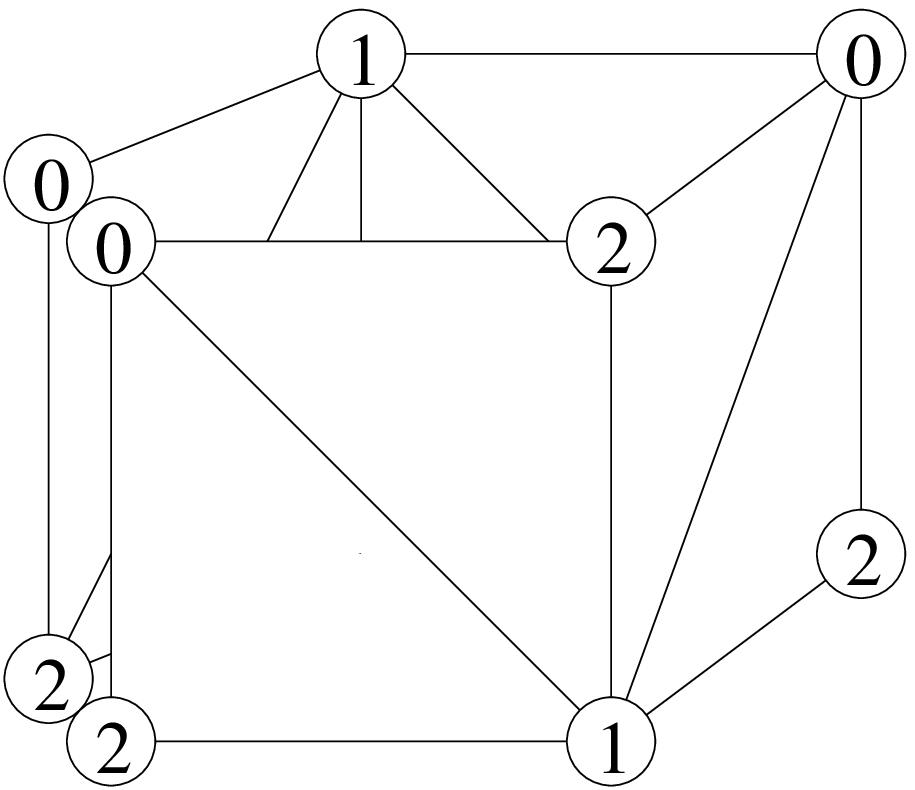}
  \raisebox{9mm}{$\rightsquigarrow$}
  \raisebox{2mm}{\includegraphics[height=16mm]{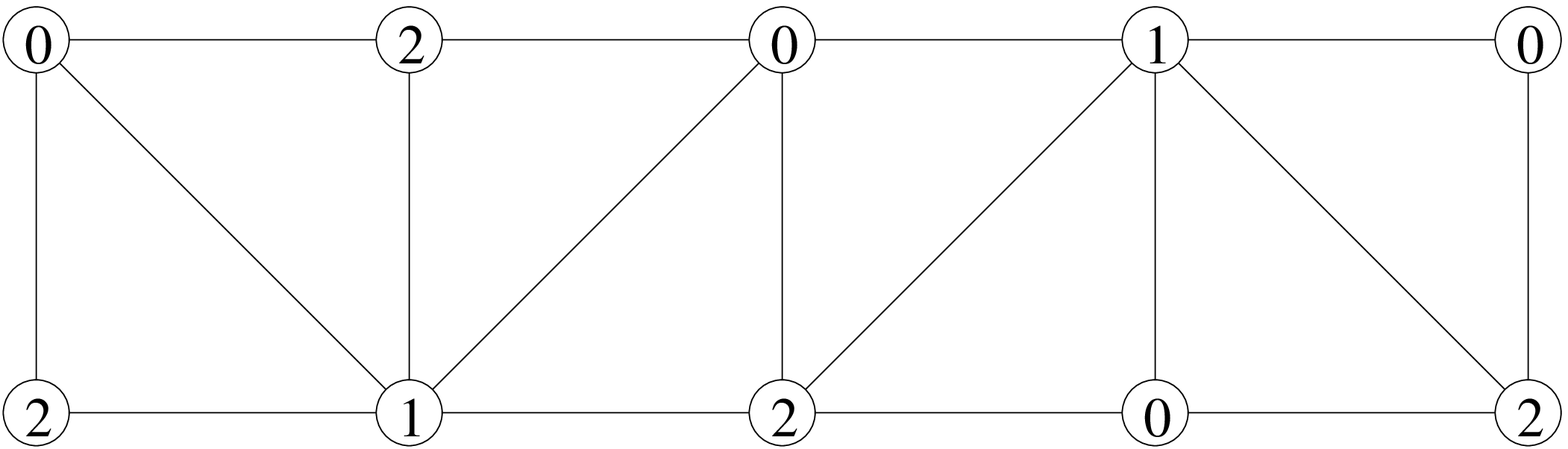}}
  \\[.5\baselineskip]
  \raisebox{9mm}{The lower chimney:}
  \includegraphics[width=25mm]{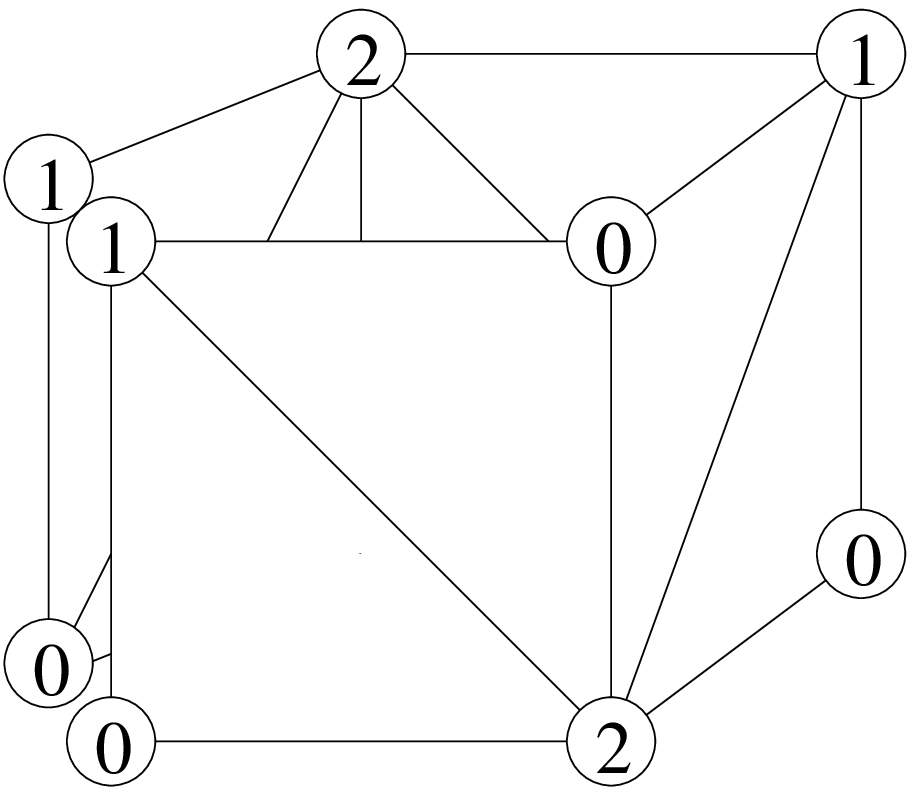}
  \raisebox{9mm}{$\rightsquigarrow$}
  \raisebox{2mm}{\includegraphics[height=16mm]{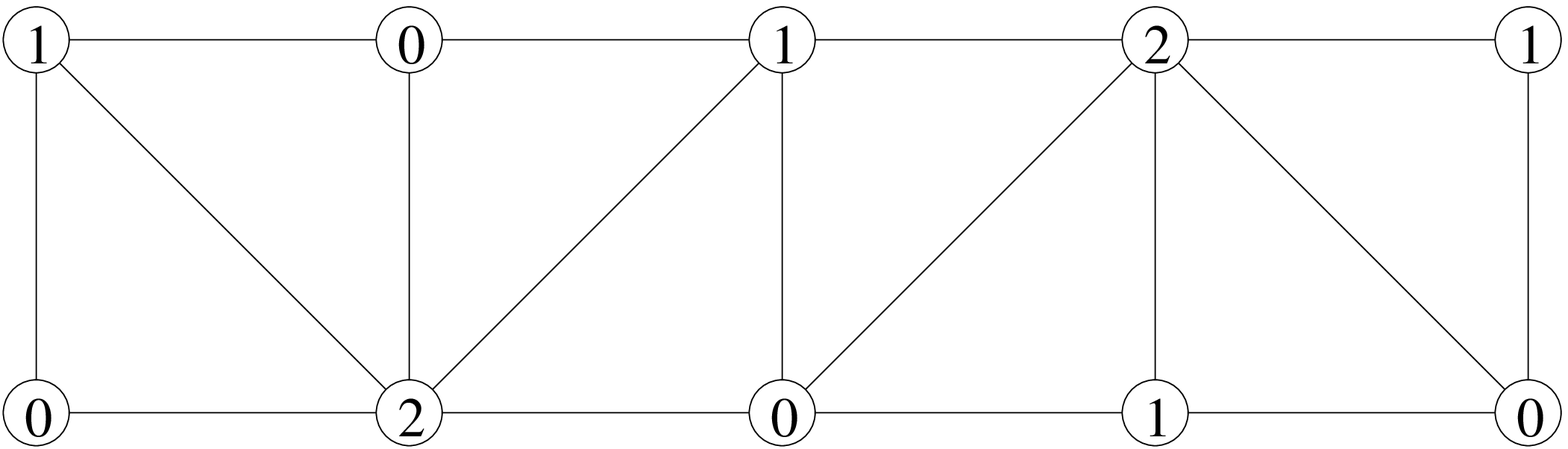}}
  \Fig{fig:chimneys}
\end{center}
We now perform a sequence of stellar subdivisions on faces of $Q\pol$, 
as follows:
\\
$\bullet$ Perform stellar subdivisions on vertical edges in the
chimneys, namely in
the upper chimney on the two edges with labels
$\begin{smallmatrix}0\\2\end{smallmatrix}$, and in
the lower chimney on the two edges with labels
$\begin{smallmatrix}1\\0\end{smallmatrix}$
(cf.~Figure~\ref{fig:pull}).
\medskip

\begin{center}
  \begin{tabular}{c c}
    The upper chimney & The lower chimney
    \\[2mm]
    \includegraphics[height=16mm]{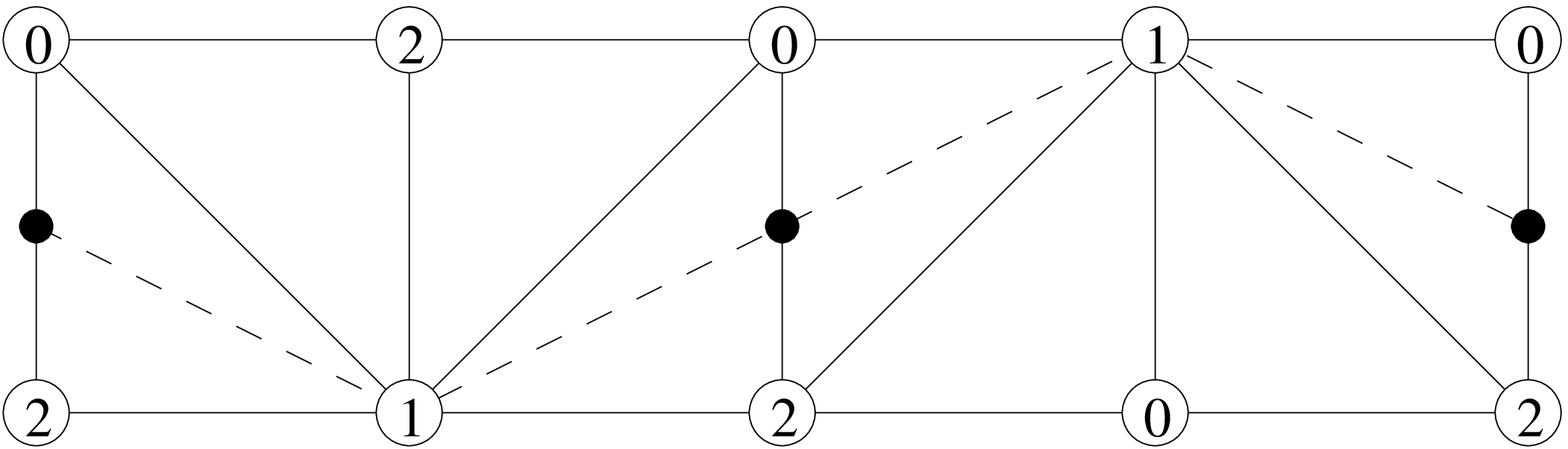} &
    \includegraphics[height=16mm]{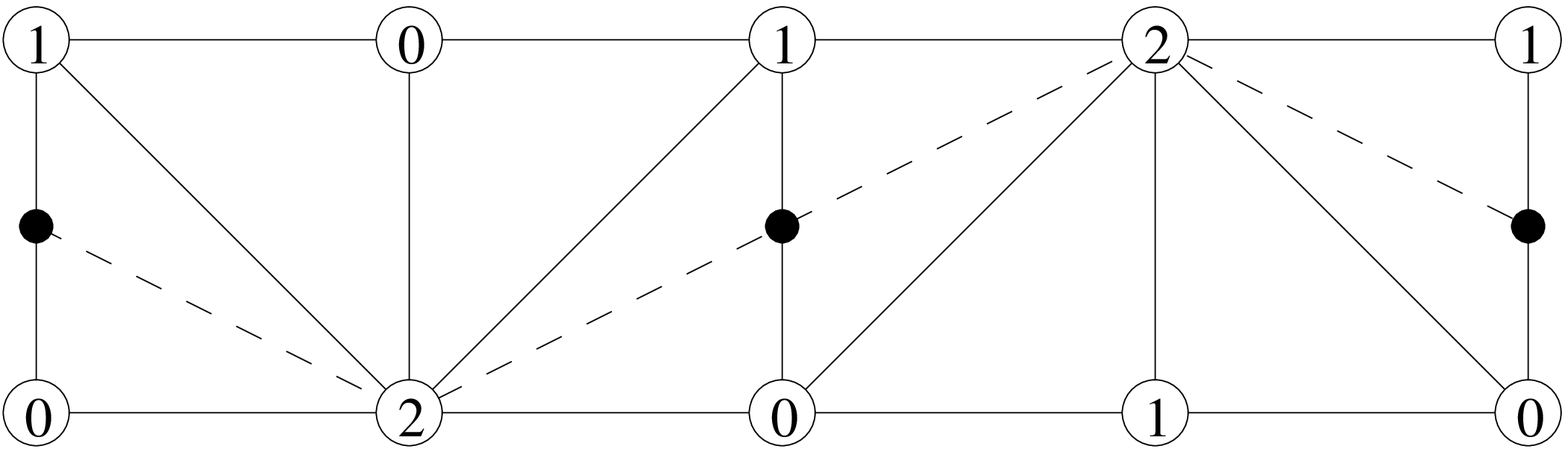}
  \end{tabular}
  \Fig{fig:pull}
\end{center}
The resulting new edges in the chimneys (four each; dashed edges in
Figure~\ref{fig:pull}) are interpreted
as marking the boundary between `inside' and `outside' in the
chimneys. No color is given to these four new vertices.

With the dashed separation edges in place,
we now have achieved the situation that every triangle
in the core is adjacent to exactly two of 
``outside'', ``upstairs'' and ``downstairs''.
Furthermore, if a triangle is adjacent to ``outside''
then it has a $0$-colored vertex, for ``upstairs'' it has
a $1$-colored vertex, and for ``downstairs'' a $2$-colored vertex.

\noindent
$\bullet$ Now perform stellar subdivisions on simplices outside $\B$ all
whose faces belong to $\B$ (e.~g., the diagonals of the cubes). Thus
$\B$ becomes an induced subcomplex.
\\
$\bullet$ Whenever there are vertices $v_1,v_2\in{\B}$ of the same
color, and a vertex $w \not\in \B$ such that both $\{w,v_1\}$ and
$\{w,v_2\}$ are edges (e.~g. $v_1,v_2$ on the boundary of the pile and
$w$ the cone vertex), perform a stellar subdivision on one of these
edges. 

The result of all these operations is now called $P\pol$. 
Thus the polytope $P$, a counterexample to Perles' original conjecture,
is obtained as its polar dual. It remains to name the
vertex star pieces that we want to attach to $\B$ in order to get a
Perles-contradicting subcomplex. For that purpose we use the vertex
labels, and proceed according to the scheme
\medskip

\begin{center}
  \begin{tabular}{c c l}\hline
    label && partial vertex star \\ \hline
    0 & $\longleftrightarrow$ & outside \\
    1 & $\longleftrightarrow$ & upstairs room\\
    2 & $\longleftrightarrow$ & downstairs room\\ \hline
  \end{tabular}
\end{center}
\medskip

If a vertex with label 2 does not touch the downstairs
room, then no partial vertex star is assigned, etc. The union $\Gamma$ 
of these partial vertex stars is a non-separating subcomplex, all
whose tetrahedra have exactly one free triangle opposite to `their'
vertex of~$\B$.\hfill\qed

\section{Connectivity \label{sec:connect}}

Finally, we verify that our counterexample does not even weakly
satisfy Perles' conjecture, that is, $\Gamma$ is dually
$3$-connected.
Our arguments will be general enough to work as well for
other counterexamples of the same type (where the core
is essentially a manifold).
Of course, for specific counterexamples, such as
the one to be presented in \url{www.eg-models.de}, $3$-connectedness
could as well be verified by explicit computer calculation.

Remove two tetrahedra $\bar{\sigma}_1,\bar{\sigma}_2$ from
$\Gamma$. We want to join any two remaining tetrahedra
$\sigma_1,\sigma_2$ by a dual path, avoiding
$\bar{\sigma}_1,\bar{\sigma}_2$. This path is constructed in three
steps:
\begin{enumerate}
\item Join $\sigma_1,\sigma_2$ to tetrahedra which have a
  full-dimensional intersection with $\core(\Gamma)$.
\item Join these full-dimensional intersections by a dual path in
  $\core(\Gamma)$.
\item Lift this path to a path in $\Gamma$.
\end{enumerate}
Step~1 follows from a simple fact about graphs of $3$-polytopes:
\begin{lemma} \label{lemma:fullIntersect}
  Every tetrahedron $\sigma \in \Gamma$ can be joined to some tetrahedron with
  full-dimensional intersection with $\core(\Gamma)$, avoiding
  $\bar{\sigma}_1,\bar{\sigma}_2$.
\end{lemma}
\begin{proof}
  The dual graph of the
  star of $v(\sigma)$ in $\Delta$ corresponds to a
  (simple) $3$-dimensional polytope, and is therefore $3$-connected. 
  Any dual path form $\sigma$ to a tetrahedron in a different component
  has to pass through $\core(\Gamma)$.
\end{proof}
For Step~2, we need some preparation. In view of Step~3, 
we consider a graph $\auxGraph$ defined as follows: its nodes
are the two-faces of $\core(\Gamma)$, and its edges correspond to
pairs $(\tau_1,\tau_2)$ of two-faces which share a one-face $\varrho$,
and which can be joined by a dual path in $\st(\varrho;\Gamma)$.
\begin{center}
  \input{vm6a.pstex_t}
  \Fig{fig:auxGraph}
\end{center}
We are interested in
the subgraphs of $\auxGraph$ that are induced by those nodes
(triangles) that lie in the star of a vertex of $\core(\Gamma)$.
\begin{lemma}
  Let $\varrho$ be a one-face of the triangle $\tau \in \core(\Gamma)$.
  Then there is at least one more triangle $\tau' \in \core(\Gamma)$
  which contains $\varrho$, and is joined to $\tau$ by an edge in
  $\auxGraph$.
\end{lemma}
\begin{proof}
  $\tau$ is covered from both sides with tetrahedra $\sigma_\uparrow$ and
  $\sigma_\downarrow\in\core(\Gamma)$. The vertices $v(\sigma_\uparrow) \neq
  v(\sigma_\downarrow)$ are both vertices of $\tau$, so that every
  edge of $\tau$ contains at least one of them.
\end{proof}
\begin{cor}\label{cor:localConnect}
  For every vertex $v\in \core(\Gamma)$, the subgraph of
  $\auxGraph$ induced by the triangles incident to $v$
  is $2$-connected.
\end{cor}
\begin{proof}
  If $v$ is a regular vertex (all incident edges
  have degree $2$ in $\core(\Gamma)$), then the corresponding local subgraph of
  $\auxGraph$ is a cycle: It is the dual graph of the star
  of~$v$ in $\core(\Gamma)$. 

  If $v$ is a singular vertex (two incident edges
  have degree $3$ in $\core(\Gamma)$), then the subgraph of
  $\auxGraph$ is a union of three paths between two nodes.
\end{proof}

Now, to verify $3$-connectivity of $\auxGraph$,
we use the following simple but very useful observation
of Naatz: 
\begin{thm}[{\cite[Theorem~3.3]{N}}]
  A graph $G$ on at least $k+1$ vertices is $k$-connected if and only 
  if for every two vertices $v$ and $w$ at distance $2$ there are at 
  least $k$ independent $v$-$w$-paths in $G$.
\end{thm}
\begin{prop}\label{prop:auxConnect}
  $\auxGraph$ is $3$-connected.
\end{prop}
\begin{proof}
  Remove two nodes/triangles $\bar{\tau}_1,\bar{\tau}_2$ from
  $\auxGraph$, and denote $\tau_0, \tau_\infty$ the two nodes which we 
  want to join by a $\bar{\tau}_1,\bar{\tau}_2$-avoiding path. In
  $\core(\Gamma)$, the intersection $\bar{\tau}_1 \cap \bar{\tau}_2$
  contains at most two vertices. 
  Now the $1$-skeleton of $\core(\Gamma)$ is $3$-connected
  (e.~g.\ by Naatz' lemma). Thus there is a 
  $\bar{\tau}_1 \cap\bar{\tau}_2$-avoiding vertex-edge-path 
  in $\core(\Gamma)$ from a 
  vertex of $\tau_0$ to some vertex of $\tau_\infty$. That is, there
  are vertices $v_0, \ldots, v_k$ of $\core(\Gamma)$ such that 
  $v_0 \in \tau_0$, $v_k \in \tau_\infty$, and the $[v_i,v_{i-1}]$'s are
  edges of $\core(\Gamma)$.
  \begin{center}
    \input{vm2a.pstex_t}
    \Fig{proof1}
  \end{center}
  Choose triangles $\tau_i \in \core(\Gamma)$ adjacent to the edges
  $[v_{i-1},v_i]$.
  \begin{center}
    \input{vm4a.pstex_t}
    \Fig{proof2}
  \end{center}
  By Corollary~\ref{cor:localConnect}, consecutive triangles $\tau_i$
  and $\tau_{i+1}$ can be joined within the star of $v_i$, avoiding
  $\bar{\tau}_1,\bar{\tau}_2$.
\end{proof}
Finally, we get Step~3 for free:
\begin{cor}\label{cor:connect}
  $\Gamma$ is dually $3$-connected.
\end{cor}
\begin{proof}
  Because $\core(\Gamma)$ is induced, $\bar{\sigma}_i \cap \core(\Gamma)$ is 
  a face of $\bar{\sigma}_i$. If $\bar{\sigma}_i \cap \core(\Gamma)$ is a
  triangle, remove the corresponding node from $\auxGraph$. 
  If $\bar{\sigma}_i \cap \core(\Gamma)$ is an edge, 
  remove the edge from $\auxGraph$, which is cut 
  by~$\bar{\sigma}_i$. By Lemma \ref{lemma:fullIntersect}, we can
  assume that we want to join two tetrahedra $\sigma_1, \sigma_2 \in\Gamma$ 
  which intersect $\core(\Gamma)$ in triangles $\tau_1$ and
  $\tau_2$ respectively. By Lemma~\ref{prop:auxConnect} we can find a path from
  $\tau_1$ to $\tau_2$ in the remaining graph. This path can be lifted
  to a dual path which joins $\sigma_1$ and $\sigma_2$, and which
  avoids $\bar{\sigma}_1, \bar{\sigma}_2$.
\end{proof}

\section{Remarks}

It seems reasonable to ask whether
every induced, $3$-connected, $3$-regular, non-separating and
planar (!) subgraph of the graph of a simple $4$-polytope
is the graph of a facet. 

Does the $120$-cell satisfy Perles' conjecture?

\subsection*{Acknowledgements}
\quad\sloppy
Thanks, in particular, to
Tom Braden,
W\l odek Kuperberg,
Carsten Lange,
Mark de Longueville,
Carsten Schultz,
Raik Stolletz and 
Elmar Vogt
for crucial discussions.

\newcommand\reference[3]{{\sc #1\ \it #2\ \rm#3}}

\vspace{\baselineskip}

\begin{thebibliography}{88}
\itemsep=0pt

\bibitem{A}
\reference
{H. Achatz \& P. Kleinschmidt:}
{Reconstructing a simple polytope from its graph,}
{in: ``Polytopes -- Combinatorics and Computation''
(G.~Kalai, G.M. Ziegler, eds.), {\sl DMV Seminar} {\bf 29}, 
Birkh\"auser, Basel 2000, pp.~155-165.}

\bibitem{Bing}
\reference
{R~H~Bing:} 
{Some aspects of the topology of $3$-manifolds related to the Poincar\'e
conjecture,}
{in: {\sl Lectures on Modern Mathematics, Vol.~II} (T.L.~Saaty, ed.),
Wiley 1964, 93-128.}

\bibitem{Bj}
\reference
{A.~Bj\"orner:}
{Shellable and Cohen-Macaulay partially ordered sets,} 
{{\sl Transactions Amer.\ Math.\ Soc.} {\bf 260} (1980), 159-183.}

\bibitem{BM}
\reference
{R. Blind \& P. Mani-Levitska:}  
{On puzzles and polytope isomorphisms,} 
{{\sl Aequationes Math.} {\bf 34} (1987), 287-297.}

\bibitem{JLP}
\reference
{G. Bohus, W. Jockusch, C. W. Lee \& N. Prabhu:}
{On a triangulation of the $3$-ball and the solid torus,}
{{\sl Discrete Math.} {\bf 187} (1998), 259-264.}

\bibitem{Borsuk}
\reference
{K. Borsuk:}
{\"Uber das Ph\"anomen der Unzerlegbarkeit in der Polyedertopologie,}
{{\sl Commentarii Math.\ Helv.} {\bf 8} (1935), 142-148.}

\bibitem{GJ1}
\reference
{E. Gawrilow \& M. Joswig:}
{{\tt polymake}{\rm:} A framework for analyzing convex polytopes,}
{in: ``Polytopes -- Combinatorics and Computation''
(G.~Kalai, G.M. Ziegler, eds.), {\sl DMV Seminar} {\bf 29}, 
Birkh\"auser, Basel 2000, pp.~43-73.}

\bibitem{GJ2}
\reference
{E. Gawrilow \& M. Joswig:}
{Polymake: A framework for analyzing convex polytopes, \rm 1997-2000,}
{\url{http://www.math.tu-berlin.de/diskregeom/polymake/doc/}.}

\bibitem{Gr}
\reference                 
{B.~Gr\"unbaum:}                 
{Convex Polytopes,} 
{Interscience, London 1967.}

\bibitem{Hatcher}
\reference
{A. Hatcher:}
{Algebraic Topology,}
{\url{http://www.math.cornell.edu/~hatcher/}, Cornell University, 2000; 
Cambridge University Press 2001, to appear.}

\bibitem{HK1}
\reference
{F. Holt \& V. Klee:}
{Counterexamples to the strong $d$-step conjecture for $d\ge5$,}
{{\sl Discrete Comput.\ Geometry} {\bf 19} (1998), 33-46.}

\bibitem{J}
\reference
{M. Joswig:}
{Reconstructing a non-simple polytope from its graph,}
{in: ``Polytopes -- Combinatorics and Computation''
(G.~Kalai, G.M. Ziegler, eds.), {\sl DMV Seminar} {\bf 29}, 
Birkh\"auser, Basel 2000, pp.~167-176.}

\bibitem{KK}
\reference
{M. Joswig, V. Kaibel \& F. K\"orner:}
{On the $k$-systems of a simple polytope,}
{Preprint \url{math.CO:0012204}, TU Berlin, December 2000, 8~pages;
{\sl Israel J. Math.}, to appear.}

\bibitem{K1}
\reference
{G.~Kalai:} 
{A simple way to tell a simple polytope from its graph,}
{{\sl J.~Combinatorial Theory, Ser.~A} {\bf 49} (1988), 381-383.}

\bibitem{K2}
\reference
{G. Kalai:}
{Some aspects of the combinatorial theory of convex polytopes,}
{in: ``{Polytopes: Abstract, Convex and Computational}''
(T.~Bisztriczky, P.~McMullen, and A.~Weiss, eds.),
Proc.\ NATO Advanced Study Institute, Toronto 1993, 
Kluwer Academic Publishers 1994, pp.~205-230.}


\bibitem{KlWa}
\reference
{V. Klee \&  D. W. Walkup:}
{The $d$-step conjecture for polyhedra of dimension $d<6$,}
{{\sl Acta Math.} {\bf 117} (1967), 53-78.}

\bibitem{Mu}
\reference
{J.~R.~Munkres:}                 
{Elements of Algebraic Topology,} 
{Addison Wesley 1984.}

\bibitem{N}
\reference
{Michael~Naatz:}                 
{The graph of linear extensions revisited,} 
{{\sl SIAM J.~Disc.~Math.} {\bf 13} (2000), 354-369.}

\bibitem{Perles}
\reference
{M. Perles:}
{Results and problems on reconstruction of polytopes,}
{Jerusalem 1970, unpublished \cite{K1}; problems posed in 
Oberwolfach 1984 and 1986 \cite{BM}.}

\bibitem{St}
\reference
{R. Stolletz:}
{Facetten-Graphen einfacher Polytope,}
{Diplomarbeit, TU Berlin 1998, 77~pages.}

\bibitem{Zeeman}
\reference
{E. C. Zeeman:}
{On the dunce hat,}
{{\sl Topology} {\bf 2} (1963), 341-358.}

\bibitem{Z1}
\reference
{G.~M.~Ziegler:}
{Lectures on Polytopes,}
{{\sl Graduate Texts in Mathematics} {\bf 152}, 
Springer-Verlag New York Berlin Heidelberg 1995;
Revised edition, 1998.}

\bibitem{Z2}
\reference
{G.~M.~Ziegler:}
{Cyclic polytopes and their duals,}
{in preparation.}

\end{thebibliography}
\end{document}